\documentclass[12pt]{amsart}

\def\squarebox#1{\hbox to #1{\hfill\vbox to #1{\vfill}}} 
\newcommand{\stopthm}{\hfill\hfill\vbox{\hrule\hbox{\vrule\squarebox 
                 {.667em}\vrule}\hrule}\smallskip} 

\newcommand{\cP}{{\mathcal P}} 
\newcommand{\pa}{\partial}
\newcommand{\Up}{\Upsilon}

\newcommand{\E}{{\mathcal E}}
\newcommand{\up}{\Upsilon}
\newcommand{\ep}{\epsilon}

\newcommand{\R}{{\mathbb R}}
\newcommand{\N}{{\mathbb N}}
\newcommand{\Hb}{{\mathbb H}}

\newcommand{\C}{{\mathbb C}}
\newcommand{\CI}{{C^{\infty}}}

\newcommand{\vol}{\operatorname{vol}}

\newcommand{\Ric}{\operatorname{Ric}}
\newcommand{\Res}{\operatorname{Res}}

\newcommand{\tr}{\operatorname{tr}}
\newcommand{\Xint}{\buildrel \circ\over X}

\renewcommand{\Re}{\mathop{\rm Re}\nolimits}

\theoremstyle{plain}
\newtheorem{theorem}{Theorem}[section]

\newtheorem{proposition}[theorem]{Proposition}

\theoremstyle{definition}

\theoremstyle{remark}

\numberwithin{equation}{section}

\title{$Q$-Curvature and Poincar\'e Metrics}
\author{Charles Fefferman}
\address{Department of Mathematics, Princeton University\\
Princeton, NJ 08544}
\email{cf@math.princeton.edu}

\author{C. Robin Graham}
\address{Department of Mathematics, University of Washington,
Box 354350\\
Seattle, WA 98195}
\email{robin@math.washington.edu}

\begin{document}
\maketitle

\thispagestyle{empty}

\renewcommand{\thefootnote}{}
\footnotetext{The work of the first author was partially supported by NSF
grant DMS-0070692.}

\section{Introduction}\label{intro}
This article presents a new definition of Branson's $Q$-curvature in
even dimensional conformal geometry.  The $Q$-curvature is a
generalization of the scalar curvature in dimension 2:  it satisfies an
analogous transformation law under conformal rescalings of the metric
and on conformally flat manifolds its integral is a multiple of the Euler
characteristic.  Our approach is motivated by the recent work 
\cite{GZ}; we derive the $Q$-curvature as a coefficient in 
the asymptotic expansion of the formal solution of a boundary problem at
infinity for the Laplacian in the Poincar\'e metric associated to the
conformal structure.  This gives an easy proof of the result
of \cite{GZ} that the log coefficient in the volume expansion of a
Poincar\'e metric is a multiple of the integral of the $Q$-curvature, 
and leads to a definition of a non-local version of the $Q$-curvature in
odd dimensions.

The $Q$-curvature is intimately connected with a family of conformally
invariant differential operators generalizing the conformal Laplacian
$\Delta +\frac{n-2}{4(n-1)}R$, for which the scalar curvature arises as the 
zeroth order term.  (Our sign convention is such that $\Delta$ is a
positive operator.)  The next operator in the family was
discovered by Paneitz \cite{Pa} and has the same principal part as 
$\Delta^2$.  Branson and \O rsted \cite{BO} observed that the zeroth order
term of Paneitz' operator gives rise to the quantity 
$$
Q = (\Delta R + R^2 - 3 | {\rm{Ric}}|^2)/6
$$
in 4 dimensions with a 
conformal transformation law similar to that of scalar curvature in
dimension 2.  This $Q$-curvature in dimension 4 
has been the focus of tremendous activity in recent years 
leading to great advances in our understanding of 4 dimensional conformal
geometry; see \cite{CY} for a survey of some of this work.  
The full family of ``conformally
invariant powers of the Laplacian'' was derived in \cite{GJMS}, and 
Branson \cite{B} formulated the definition of the $Q$-curvature in general 
even dimensions using these operators.  However, this general definition
involves an 
analytic continuation in the dimension and the higher-dimensional
$Q$-curvature has remained a rather mysterious object.

The work \cite{GZ} shows that the conformally invariant powers of the
Laplacian and the $Q$-curvature arise naturally in
scattering theory for Poincar\'e metrics associated to the conformal
structure.  This connection can be thought of as providing a 
different definition of the $Q$-curvature, in which the analytic
continuation in the dimension is replaced by an analytic continuation in a
spectral parameter with the dimension fixed.  Conceptually this
seems an improvement, but the definition is still indirect.  Our
contribution is to provide a direct definition in the context of the
Poincar\'e metric. 

The Poincar\'e metrics involved are 
asymptotically Einstein metrics whose conformal infinity is the given
conformal structure.  The model is hyperbolic space
$\Hb^{n+1} = \{(x,y)\in \R_+\times \R^n\}$, with metric
$g=x^{-2}(dx^2 + \sum dy_i^2)$.  On $\Hb^{n+1}$, the function $U=\log x$
satisfies $\Delta_g U =n$.  In \S 3, we show that if $g$ is a formal
Poincar\'e metric associated to the conformal structure, the equation
$\Delta_g U=n$ has a formal solution asymptotic to $\log x +O(x)$, where 
$x$ is a defining function for the boundary at infinity associated to a
metric in the conformal class. 
The $Q$-curvature is a multiple of the coefficient of $x^n\log x$ in the
expansion of $U$.  The conformal transformation law for $Q$ follows
immediately from the characterization given in \cite{GZ} of the conformally 
invariant $(n/2)^{\text{th}}$ power of the Laplacian $P_{n/2}$ as the
coefficient 
of $x^n \log x$ in the expansion of a formal solution of $\Delta_g u=0$
with prescribed boundary value.  The result that the log term
coefficient in the volume expansion is the integral of $Q$ follows by an
application of Green's identity.

In \S 4, we study the equation $\Delta_g U=n$ globally.  We show that 
in any dimension there is a unique global
solution with the prescribed leading asymptotics, given by the formula 
$U=-\frac{d}{ds}\cP (s)1|_{s=n}$, where $\cP(s)$ is the family of Poisson
operators for $\Delta_g - s(n-s)$.  If $n$ is
odd, there is no $x^n \log x$ term in the expansion of $U$, but the $x^n$
coefficient provides a non-local analogue of $Q$.  This odd-dimensional
version satisfies a conformal transformation law just like the one in even
dimensions, in which $P_{n/2}$ is replaced by a pseudodifferential operator
acting invariantly on conformal densities with principal part the same
as that of $\Delta^{n/2}$; namely a multiple of $S(n)$, where $S(s)$ is
the scattering matrix for $\Delta_g - s(n-s)$.   The same Green's identity
argument as in even dimensions shows that in odd dimensions, the
renormalized volume itself is a multiple of the integral of the
odd-dimensional $Q$.  One important distinction between even and odd
dimensions, though, is that in odd dimensions both $Q$ and the operator
$S(s)$ are not determined solely by the conformal structure on the boundary
at infinity, but also depend on the extension of the formal
Poincar\'e metric to a metric on the interior.  This dependence 
could be a limitation for their use in the study of conformal geometry. 

We include in \S 3 a different proof of another result of \cite{GZ}: 
the self-adjointness of the conformal powers of the Laplacian of
\cite{GJMS}.  In \cite{GZ} this was shown to follow by analytic continuation
from the self-adjointness of the scattering matrix, of which the conformal
powers of the Laplacian arise as residues.  We show that a renormalized
integration by parts identity giving self-adjointness of the scattering
matrix can be formulated 
directly at the limiting value and gives a proof based purely on
formal asymptotics.

\section{Background}\label{backg}
In this section we discuss background material about Poincar\'e
metrics, conformally invariant powers of the Laplacian, $Q$-curvature, and 
volume renormalization which we will need later.

Let $M$ be a manifold of dimension $n>2$ with a conformal structure
$[h]$.  Any two metrics 
$h,\hat{h}\in [h]$ are conformally related:  $\hat{h}=e^{2\Up}h$ for some 
$\Upsilon \in \CI ( M )$.  The conformal structure is equivalent to the
specification of the {\em metric bundle} of $[h]$.  This is the ray bundle 
${\mathcal G}\subset S^2T^*M$ of multiples of the metric:  the fiber of
${\mathcal G}$ over $p\in M$ is 
$\{t^2h(p):t>0\}$. The space of conformal densities of weight $w\in \C$ is 
\[ \E ( w ) = \CI( M ; {\mathcal G}^{-\frac{w}2} ) \,, \]
where by abuse of notation we have denoted by ${\mathcal G}$ also the line
bundle associated to the ray bundle defined above.  
A choice of representative $h$ for the conformal structure induces an
identification $\E (w)\simeq \CI(M)$; if $\hat{h}=e^{2\up}h$ then the 
corresponding elements of $\CI (M)$ transform by $\hat{f}=e^{w\up}f$. 

Our paper \cite{FG} introduced two formal metrics in higher dimensions
which can be associated to the conformal structure $(M,[h])$.  The 
relevant one here is the Poincar\'e metric.  
If $X$ is an $n+1$-manifold
with $ \partial X = M $ and $ x $ is a defining function of 
$ \partial X =M $ in $ X$:
\[ x|_{\Xint}>0 \,, \ \ x |_{\partial X } = 0 \,, \ \ 
dx |_{ \partial X } \neq 0 \,, \]
we say that a metric $ g $ on $\Xint$ is 
conformally compact with conformal infinity $ [h] $ if
$$
 g = \frac {\overline g }{ x^2 } \,, \ \ {\overline g}|_{T\partial X} \in
[h] \,,
$$
where ${\overline g}$ is a smooth metric on $X$. 
A conformally compact metric is said to be asymptotically
hyperbolic if its sectional curvatures approach $-1$ at $\pa X$; this is
equivalent to $|dx|_{\bar g}=1$ on $\pa X$.  
It was shown in \cite{GL} that if $g$ is an asymptotically hyperbolic
metric on $X$, then a choice of metric $h$ in the conformal class on $M$
uniquely determines a defining function $x$ near $\pa X$ and an
identification of a neighborhood of $\pa X$ with $M \times
[0,\epsilon)$ such that $g$ has the normal form 
\begin{equation}\label{normalform}
g=x^{-2}(h_x + dx^2),\qquad h_0=h,
\end{equation}
where $h_x$ is a 1-parameter family of metrics on $M$.  We will say
that $x$ is the defining function associated to $h$.

The Poincar\'e metric of $(M,[h])$ is conformally compact and is an
asymptotic solution of the Einstein equation 
$\Ric(g)=-ng$.  It is easy to see that any such metric is asymptotically
hyperbolic.  In order to formally solve this equation, we identify a
neighborhood of $M$ in $X$ with  
$M\times [0,\ep)$ and consider metrics of the form
\eqref{normalform}.   
The equation $\Ric(g)+ng=0$ can be calculated directly in terms of 
$h_x$ and the formal asymptotics of solutions studied; see
\cite{G}.  If $n$ is odd, then there is a unique formal smooth solution
$h_x$ to 
$
\Ric(g)+ng= O (x^{\infty})
$
which is even in $x$.  If $n$ is even,  
the condition 
$
\Ric(g)+ng= O (x^{n-2})
$
uniquely determines 
$h_x \mod  O (x^n)$, and this $h_x$ is even in $x$ (mod $ O (x^n)$).
Although in 
general smooth solutions do not exist to higher orders for $n$ even,
the condition 
$
\tr_g(\Ric(g)+ng)= O (x^{n+2})
$
can be satisfied and uniquely 
determines the $h$-trace of the $x^n$ coefficient in $h_x$.
The indicated Taylor coefficients of $h_x$ are determined inductively from 
the equation and are given by polynomial formulae in terms of $h$, its
inverse, and its curvature tensor and covariant derivatives thereof.

Since any asymptotically hyperbolic metric can be put uniquely in the form
\eqref{normalform} upon choosing $h$, it follows that the equivalence class
of the solution $g$ up to 
diffeomorphism and up to terms vanishing to the indicated orders is
uniquely determined by the conformal structure.  This
equivalence class is
called the formal Poincar\'e metric associated to $[h]$.  When $n$ is even,
the higher order terms in $h_x$ are not determined and will not play a role
in this paper.  Similarly, it is not crucial that the metric be smooth to all
orders--honest Einstein metrics will have log terms in their expansions,
but these also occur to orders that do not affect the considerations here. 
When $n$ is odd, a number of our results hold without the
assumption that 
$h_x$ is even beyond the order forced by the Einstein equation.
However, for simplicity in statements below, we shall restrict attention 
to $h_x$ which are smooth and even in $x$ to all orders.  
If $X$ is a manifold with $\partial X=M$, we shall call a Poincar\'e metric
associated to $[h]$ any metric on $X$
whose restriction to a collar neighborhood of $M$ is such a formal
Poincar\'e metric.

The conformally invariant powers of the Laplacian are a family $P_k$, 
$k\in \N$ and $k\leq n/2$ if $n$ is even, of 
natural differential operators on $n$-dimensional Riemannian manifolds;
$P_k$ has the same principal part as $\Delta^k$ 
and equals $\Delta^k$ if $h$ is flat.
$P_k$ is conformally invariant in 
the sense that it defines an operator $P_k:\E(-n/2+k)\rightarrow
\E(-n/2-k)$ on conformal densities independent of the choice of conformal
representative, which is equivalent to the statement 
that under conformal rescaling $\hat{h}=e^{2\Up}h$, $P_k$ transforms by
$$
\widehat{P} _k =  e^{ ( - n/2 - k ) \Upsilon } 
P_k\, e^{ ( n/2 - k ) \Upsilon } \,. 
$$
The operators $P_k$ were constructed in \cite{GJMS} using the ambient
metric of \cite{FG}.  In \cite{GZ} it was 
shown that the derivation can be reformulated in terms of 
the Poincar\'e metric by consideration of the formal 
asymptotics of solutions of 
\begin{equation}\label{infiniteorder}
(\Delta_g -s(n-s))u=O(x^{\infty})\,.
\end{equation}
If $g$ is a Poincar\'e metric for $(M,[h])$ and $x,\hat{x}$ the
are the defining functions associated to the representative metrics
$h,\hat{h}$, then $\hat{x}=e^{\Up}x +O(x^2)$.  It follows that if $f\in
\CI(M)$ represents 
a section of $\E (w)$, then it is a conformally invariant statement to
require that a function $u$ on $X$ be asymptotic to $x^{-w}f$.  
The operator $P_k$ is derived in \cite{GZ} as the obstruction to the
existence of a smooth solution
of \eqref{infiniteorder} for $s=n/2+k$ with $u$ having prescribed leading
term $x^{n/2-k}f$:

\begin{proposition}(\cite[Propositions 4.2, 4.3]{GZ})
\label{conformalexpansions}
Let $g$ be a Poincar\'e metric associated to $(M,[h])$ and let  
$f\in \CI (M)$.  Let $k\in \N$ with $k\leq n/2$ if $n$ is even. 
There is a solution of \eqref{infiniteorder} for $s=n/2+k$
of the form  
\begin{equation}
\label{logform}
u=x^{n/2-k}(F + H x^{2k} \log x)
\end{equation}
with $F,H\in \CI(X)$ and with $F|_M=f$.
The functions $F \mod O  (x^{2k})$ and $H|_M$ are formally determined by
$h$, and 
\begin{equation}\label{logterm}
H|_M = -2c_k P_kf\,,\qquad c_k = (-1)^k[2^{2k}k!(k-1)!]^{-1}\,,
\end{equation}
where $P_k$ is the conformally invariant operator of \cite{GJMS}.
\end{proposition}

Branson's definition of the $Q$-curvature is in terms of the zeroth order
terms of the operators $P_k$.  
Note that if $n$ is even, then $P_{n/2}$ is invariantly defined from
$\E(0)=\CI(M)$  to $\E(-n)$.  Since the constant function $1\in \E(0)$ has
a smooth extension annihilated by $\Delta_g$, we have $P_{n/2}1=0$.  
Therefore $P_{n/2}$ has no constant term in dimension $n$.  However, if we
fix $k\in \N$ and denote by $P_k^n$ the operator $P_k$ in dimension $n$,
then $P_k^n$ has coefficients given by universal formulae in curvature and
its derivatives which are rational in $n$, and the zeroth order term of 
$P_k^n$ is of the form $(n/2-k)Q^n_k$ for a scalar Riemannian invariant
$Q_k^n$ with coefficients rational in $n$ and regular at $n=2k$.  
The $Q$-curvature in even dimension $n$ is then defined as $Q=Q^n_{n/2}$.   
An analytic continuation argument in the dimension using the conformal
transformation law of the operators $P_k$ (see \cite{B}) yields the
invariance property of $Q$:  if once again $\hat{h} = e^{ 2 \Upsilon } h$,
then  
\begin{equation}
\label{Qtrans}
e^{ n \Upsilon } {\widehat Q} = Q + P_{n/2}\Upsilon \,.
\end{equation}
From this and the facts that $P_{n/2}$ is self-adjoint and $P_{n/2}1=0$, it
follows that $\int_MQ\,dv_h$ is independent of the choice of $h$ in the
conformal class.  

We have learned that Gover-Peterson \cite{GP} have found a construction of
the operators $P_k$ and $Q$ using the conformal tractor calculus.

Let $(X,g)$ be a Poincar\'e metric associated to $(M,[h])$.  The
AdS/CFT correspondence in physics has given rise to interesting invariants 
which arise in the asymptotics of the volume function on $X$ -- 
see \cite{W}, \cite{HS}, \cite{GW}, \cite{G}.  To define these invariants, 
observe that if $h$ is a metric in the conformal class and $x$ is the
associated defining function, then \eqref{normalform} implies that near
$\pa X$, the volume form of $g$ takes the form
$dv_g =x^{-n-1} dv_{h_x}dx$,
where
\begin{equation}\label{volhx}
dv_{h_x}=
(1+v^{(2)}x^2+(\mbox{even powers})+\ldots)dv_h
\end{equation}
and the $v^{(2j)}$ are smooth functions on $M$
determined by $h$ (for $1\leq j\leq n/2$ if $n$ is even).  
Integration shows that 
\begin{equation}
\label{volexp}
 \begin{split}
&  \vol_g ( \{ x > \epsilon \} ) = c_0 \epsilon^{-n} + c_2 \epsilon^{-n+2}
+ \cdots + c_{n-1} \epsilon^{-1} + V + o(1) \\ 
& \ \\
& \ \ \ \ \ \ \text{for  $ n$  odd, }\\
& \ \\
&  \vol_g ( \{ x > \epsilon \} ) = c_0 \epsilon^{-n} + c_2 \epsilon^{-n+2}
+ \cdots + c_{n-2} \epsilon^{-2} 
 +  L \log ( 1 / \epsilon ) 
+ V + o(1) \\  
& \ \\
& \ \ \ \ \ \ \text{for $ n$ even. }
\end{split}
\end{equation}
It can be shown (see \cite{G}) that $V$ is independent of
the conformal representative $h$ on the boundary at infinity when $n$ is
odd, and $ L$ is independent of the conformal representative when $n$ is
even.  Graham-Zworski \cite{GZ} have recently identified $L$ as the
integral of the $Q$-curvature:  $L=2c_{n/2}\int_MQ\, dv_h$,
where $c_{n/2}$ is as in \eqref{logterm}.
We shall give another proof of this result below.
Anderson \cite{A} has identified $V$ when $n=3$.
In an appendix to \cite{PP}, Epstein shows that for conformally compact
hyperbolic manifolds,
the invariants $L$ for $n$ even and $V$ for $n$ odd are each
multiples of the Euler characteristic $\chi(X)$.
The dependence of $ V $ on the choice of $h$ for $n$ even is the
so-called holographic anomaly -- see \cite{HS},\cite{G}.  

\section{Formal Considerations}

The considerations of this section are all formal near the
boundary, so we shall take $X=M\times [0,1)$.  

\begin{theorem}\label{formal}
Suppose $n$ is even and let $(X,g)$ be a Poincar\'e metric
associated to $(M,[h])$.  Choose a representative metric $h$ with special
defining function $x$.
There is a unique solution $U \mod O(x^n)$ to
\begin{equation}\label{Uequation}
\Delta_g U = n +O(x^{n+1}\log x)
\end{equation}
of the form
\begin{equation}\label{Uexp}
U=\log x + A + Bx^n \log x +O(x^n)\,,
\end{equation}
with
$$A,B\in C^{\infty}(X)\,, \qquad A|_M=0\,. $$ 
Also, $A \mod O(x^n)$ and $B|_M$ are formally determined by
$h$. 
\end{theorem}

\begin{proof}
Write $g$ in the form \eqref{normalform}.  A straightforward calculation
shows that 
\begin{equation}\label{Laplacian}
\Delta_g = -(x\pa_x)^2 + (n-\frac{x}{2}\tr_{h_x}h'_x)x\pa_x +
x^2\Delta_{h_x}\,,
\end{equation}
where $h'_x=\pa_x h_x$.  Therefore the equation $\Delta_g U = n$ is
equivalent to
$$\Delta_g (U-\log x) = \frac{x}{2}\tr_{h_x}h'_x.$$  
In order to construct a solution to this equation, observe that 
\eqref{Laplacian} implies that if $a_j\in \CI(M)$, then
$$
\Delta_g(a_jx^j)=j(n-j)a_jx^j + O(x^{j+1}).
$$
Since $x\tr_{h_x}h'_x \in x^2\CI(X)$ and is formally determined by $h$
$\mod O(x^{n+1})$, we can take $a_0=a_1=0$ and then inductively
determine $a_j$ uniquely for $1<j<n$ so that $A=\sum_{j=2}^{n-1}a_jx^j$
solves 
$$
\Delta_g A = \frac{x}{2}\tr_{h_x}h'_x + x^n E\,, 
$$  
where $E \in \CI(X)$ and $E|_M$ is formally determined by $h$.  The $x^nE$
error cannot be removed by a term of the form $a_n x^n$, but it can by a
term of the form $B x^n\log x $ since
$$
\Delta_g(B x^n \log x)=-nBx^n + O(x^{n+1}\log x)\,.
$$
Also $B|_M$ is formally determined by $h$ since we must
take $-nB|_M = E|_M$.
\end{proof}

We observe that since $h_x$ is an even function of $x$, it
follows from the inductive construction above that $A \mod O(x^n)$ is also
even.  An argument similar to the above shows that when $n$ is odd, there
is a unique formally determined solution $U \mod O(x^{\infty})$ to 
$\Delta U =n +O(x^{\infty})$ of the form $U=\log x +A$ with $A\in \CI(X)$,
$A|_M=0$, $A \mod O(x^{\infty})$ even in $x$.  However, as we shall see in
the next section, for $n$ odd the interesting term analogous to $Bx^n\log
x$ is globally determined.  When $n$ is even, the solution $U$ of
Theorem~\ref{formal} can be continued to higher order, but is no
longer formally determined by $h$.

The function $B|_M$ in Theorem~\ref{formal} is an invariant of the
metric $h$.  We shall show in Theorem~\ref{global} that 
\begin{equation}\label{B=Q}
B|_M= -2c_{n/2}Q\,,
\end{equation}
where $Q$ is Branson's $Q$-curvature, so that the construction of
Theorem~\ref{formal} provides a characterization of $Q$.  
For now we take \eqref{B=Q} as a definition and 
show how to deduce the transformation law \eqref{Qtrans}.

Suppose we choose a different metric $\hat{h} = e^{ 2 \Upsilon } h$ in
Theorem~\ref{formal}.  Then $\hat{x}= e^{\Up}x \mod O(x^2)$ so that 
$\log \hat{x} = \log x + \Up +O(x)$.  It follows that 
$\Delta_g(\hat{U}-U) = O(x^{n+1}\log x)$ and 
$\hat{U}-U = F+ Hx^n\log x +O(x^n)$, where  
$F,H\in \CI(X)$ and $F|_M=\Up$.  By
Proposition~\ref{conformalexpansions}, this problem
solved by $\hat{U}-U$ is precisely the one
giving rise to the invariant operator $P_{n/2}$, so we deduce that 
$$
H|_M = -2c_{n/2}P_{n/2}\Up\,.
$$
However, subtracting the expansions \eqref{Uexp} for $\hat{U}$ and $U$
shows that 
$$H|_M =e^{n\Up}\hat{B}|_M - B|_M,$$ 
giving \eqref{Qtrans}.  

We next show how to use the function $U$ to give a simple proof of the
following theorem of \cite{GZ}.  Recall that $L$ denotes the coefficient of
the logarithmic term in \eqref{volexp}.

\begin{theorem}\label{L}
If $n$ is even, then
$$
L = 2 c_{n/2} \int_M Q \;dv_h \,.
$$
\end{theorem}
\begin{proof}
Let $x_0>0$ be small.  Recalling \eqref{normalform}, Green's identity gives 
$$
\int_{\ep<x<x_0}\Delta_g U \;dv_g = \ep^{1-n} \int_{x=\ep}
\pa_xU\;dv_{h_{\ep}} 
- x_0^{1-n}\int_{x=x_0}\pa_x U \;dv_{h_{x_0}} \,,
$$
which together with \eqref{Uequation} yields
$$
n\vol_g\{\ep<x<x_0\} = \ep^{1-n}\int_M \pa_xU|_{x=\ep}dv_{h_{\ep}} + O(1)\,.
$$
By definition, the coefficient of $\log \ep$ in the expansion of the left
hand side is  $-nL$.  Differentiating \eqref{Uexp} and using \eqref{volhx}
shows that the  
coefficient of $\log \ep$ on the right hand side is $n\int_M B \,dv_h$,
so the result follows from \eqref{B=Q}.
\end{proof}

The fact that the conformally invariant operators $P_k$ of \cite{GJMS} are
self-adjoint was 
proved in \cite{GZ} by analytic continuation in the spectral parameter,
using the self-adjointness on the real axis of the scattering matrix.  
One way that the scattering matrix can be proved to be self-adjoint is via
a renormalized integration by parts identity (Proposition 3.3 of \cite{GZ})
originally discussed for hyperbolic space in \cite{W}.  
In the next proposition, we show that a variant of this identity gives
self-adjointness of the $P_k$ directly.

\begin{proposition}\label{logpart}
Let $(X,g)$ be a Poincar\'e metric associated to $(M,[h])$, let
$k\in \N$ with $k\leq n/2$ if $n$ is even, and set $s=n/2+k/2$.  
Let $f_1,f_2\in \CI(M)$ and let $u_1,u_2$
denote the corresponding solutions of \eqref{infiniteorder} given by 
Proposition~\ref{conformalexpansions}.  Then for fixed small $x_0>0$,
\[
\begin{split}
 {\rm{lp} } \;
 \int_{ \ep<x<x_0}&[\langle du_1,du_2\rangle -s(n-s)u_1\overline{u_2}] 
\;dv_g \\
&=2nc_k\int_{M} f_1\overline{P_{k}f_2} \;dv_h 
= 2nc_k \int_{M} \overline{f_2}P_kf_1 \;dv_h\,,
\end{split}
\]
where $ {\rm{lp}} $ denotes the coefficient of $\log {\ep}$ in the
asymptotic expansion of the integral as $\ep \rightarrow 0$.
In particular, $P_k$ is self-adjoint.
\end{proposition}
\begin{proof}
With $g$ in the form \eqref{normalform}, Green's identity gives
\[
\begin{split}
&\int_{ \ep<x<x_0} [\langle du_1,du_2\rangle -s(n-s)u_1 \overline{u_2}] 
\;dv_g \\ 
 & =\int_{ \ep<x<x_0}u_1(\Delta  -s(n-s) )\overline{u_2} dv_g 
 +x_0^{1-n}\int_{x=x_0}u_1\pa_x\overline{u_2} dv_{h_{x_0}}
-\ep^{1-n}\int_{x=\ep}u_1\pa_x\overline{u_2} dv_{h_{\ep}}\\
&= -\ep^{1-n}\int_{x=\ep}u_1\pa_x\overline{u_2}\;dv_{h_{\ep}} + O(1)\,. 
\end{split}
\]
Substituting \eqref{logform} and \eqref{volhx} and expanding shows that the 
coefficient of $\log \ep$ in the expansion of this expression is
$$
-\int_M\left[ (n/2+k)F_1\overline{H_2} + (n/2-k)H_1\overline{F_2}\right]
\;dv_h\,.
$$
By symmetry in $u_1,u_2$, we deduce that this equals
$$
-n\int_M F_1\overline{H_2}\;dv_h = -n\int_M H_1\overline{F_2}\;dv_h \,,
$$
which gives the desired result by \eqref{logterm}.
\end{proof}

\section{Global Considerations}
In this section, let $X$ be a
compact $n+1$-manifold with $\pa X=M$.  (If there is no such $X$, one can
take instead $X=M\times [0,1]$, for which $\pa X = M \sqcup M$.)  Let $g$
be a Poincar\'e metric on $X$ with conformal infinity 
$(M,[h])$, and choose a metric $h$ in the conformal class with associated
defining function $x$.  The main theorem of this section is the following:

\begin{theorem}\label{global}
There is a unique function $U\in \CI(\Xint)$ solving
$$
\Delta_g U = n
$$
and with the asymptotics
\begin{equation}\label{Uasym}
U=\left\{
\begin{array}{ll} 
\log x + A +Bx^n\log x & \text{ for $ n $ even}\\ 
\log x + A + Bx^n &  \text{ for $ n $ odd,} 
\end{array} \right.
\end{equation}
where $A,B\in \CI(X)$ are even $\mod O(x^{\infty})$ and $A|_M=0$.
Moreover, if $n$ is even, we have
\begin{equation}\label{B=Q2}
B|_M= -2c_{n/2}Q\,.
\end{equation}
\end{theorem}

Our analysis of this global problem depends on results from \cite{GZ} on
the asymptotic Dirichlet problem at infinity for the operator 
$\Delta_g - s(n-s)$, which we briefly review.  Mazzeo and Mazzeo-Melrose 
\cite{Ma1}, \cite{MM}, \cite{Ma2} showed that the 
spectrum of $\Delta_g$ is of the form $\Sigma \cup [(n/2)^2,\infty)$, 
where $\Sigma \subset (0,(n/2)^2)$ is the finite set of $L^2$ eigenvalues.  
Using fundamental results of \cite{MM} concerning the existence and
properties of the resolvent, in \cite{GZ} (see also \cite{JS}), it was
shown that there is a meromorphic family of Poisson operators 
$$ \cP ( s): \E(s-n) \longrightarrow \CI(\Xint)$$ 
for $\Re s > n/2$, with poles only for $s$ such that $s(n-s)\in \Sigma$,
with the characterizing properties:
$$(\Delta_g -s(n-s))\cP(s)f = 0\,,$$
$$
\begin{array}{ll}
\cP(s)f = x^{n-s}F + x^s G & \mbox{if } s\notin n/2+\N \\
\cP(s)f = x^{n/2-k}F + H x^{n/2+k}\log x & \mbox{if } s=n/2+k, \,
k\in\N\,,
\end{array}
$$
for $F,G,H\in \CI(X)$ such that $F|_M=f$.  If $2s-n$ is not an odd 
integer, then  
$F,G,H \mod O(x^{\infty})$ are uniquely determined 
and are even.  If $2s-n$ is an odd integer, then 
$F,G \mod O(x^{\infty})$ may be chosen to be even and are then uniquely
determined.   

If $s=n/2+k$ for $k\in \N$ satisfying $k\leq n/2$ if $n$ is even,
then $H|_M$ is locally determined by $f$ and $h$.  In fact, by
Proposition~\ref{conformalexpansions}, in this case we have $H|_M = -2c_k
P_kf$ where $P_k$ is the conformally invariant operator on $M$.  However,
if $s\notin n/2+\N$, then $G|_M$ is globally determined by $f$ and $g$ and
defines the scattering matrix:  
\begin{equation}\label{scat}
S(s)f = G|_M\,.
\end{equation}
It can be shown 
(\cite{Me}, \cite{JS}, \cite{GZ}) that 
$$
  S (s ) : \E(s-n) \longrightarrow  \E(-s)   
$$
is a meromorphic family of pseudodifferential operators on $M$ for 
$\Re s >n/2$, of order $2\Re s -n$, invariantly
defined on densities, with poles only for $s(n-s)\in
\Sigma$ and $s\in n/2+\N$, and with principal symbol
\begin{equation}\label{symbol}
\sigma ( S ( s) ) = 2^{ n - 2s } \frac{ \Gamma( n/2 - s) }{ \Gamma ( 
s - n /2 ) }\sigma ( \Delta_h ^{s - n/2 } ) \,.
\end{equation}
The full symbol of $S(s)$ is determined by $g \mod O(x^{\infty})$, and
therefore by $h$ if $n$ is odd (and to the appropriate order if $n$ is
even).   
Theorem 1 of \cite{GZ} shows that if $k\in \N$ and $k\leq n/2$ for $n$
even, and $(n/2+k)(n/2-k)\notin \Sigma$, then the pole of $S(s)$ at
$s=n/2+k$ is simple and 
\begin{equation}\label{residue}
\Res_{s=n/2+k}S(s)=-c_kP_k.
\end{equation}
Here we have used the metric $h$ to trivialize the 
density bundles and so regard all $S(s)$ as operators on $\CI(M)$.

Observe that $\cP(s)$ is always regular at $s=n$, $\cP(n)$ acts invariantly
on $\CI(M)=\E(0)$, and we have $\cP(n)1=1$.  
In particular, from the expansion we obtain $P_{n/2}1=0$ for $n$ even and
$S(n)1 = 0$ for $n$ odd.  If $n$ is even, \eqref{residue} then shows 
that $S(s)1$ extends holomorphically across $s=n$, and if we set 
$S(n)1=\lim_{s\rightarrow n}S(s)1$, then Theorem 2 of \cite{GZ} states that 
\begin{equation}\label{Sn1}
S(n)1 = c_{n/2} Q.
\end{equation}

Still using $h$ to trivialize the density bundles, consider the expansion 
\begin{equation}\label{exp1}
\cP(s)1 = x^{n-s}F_s + x^s G_s
\end{equation}
for $s$ near $n$, with $F_s, G_s \mod O(x^{\infty})$ even as described
above.
Of course, such an expansion holds also for $s=n$, and for $n$ odd it is
clear that $F_n=1$, $G_n=0 \mod O(x^{\infty})$.  
However, if $n$ is even, the expansion alone does not uniquely determine
$F_n, G_n \mod O(x^{\infty})$.
The construction of the Poisson operators in \cite{GZ} shows that 
$F_s, G_s\in \CI(X)$ may be chosen to depend holomorphically on $s$ across
$s=n$, which fixes this ambiguity in $F_n, G_n \mod O(x^{\infty})$. 
The definition \eqref{scat} of the scattering matrix together with
\eqref{Sn1} give $G_n|_M = c_{n/2}Q$, so the condition $\cP(n)1=1$ forces  
$F_n = 1 - c_{n/2} Q x^n + O(x^{n+1})$. 
(Actually, the proof of \eqref{Sn1} in \cite{GZ} is by the reverse of this
logic:  the $x^n$ coefficient of $F_n$ is evaluated first by taking the
limit of the formal solutions $F_s \mod O(x^{\infty})$. By arguing along
these lines, it is possible to give a purely formal proof of \eqref{B=Q2}.) 

\medskip\noindent
{\em Proof of Theorem~\ref{global}.}
Consider the holomorphic family of functions $\cP(s)1\in \CI(\Xint)$ for
$s$ near $n$.  We have 
$[\Delta_g -s(n-s)]\cP(s)1 = 0$. 
Differentiation with respect to $s$ at $s=n$ gives
$\Delta_g U = n$, where 
$$
U = -\frac{d}{ds}\cP (s)1|_{s=n}\,.
$$
Differentiation of \eqref{exp1} yields
$$
U = F_n\log x  - F_n' - G_n x^n\log x  - x^n G_n'\,,
$$
where $'=d/ds$.
Substituting the values derived above (and recalling that $F_s|_M=1$ so
that $F_n'|_M=0$) gives the stated asymptotics of $U$.

Uniqueness follows from the facts that the asymptotics of $U$ up to the 
$x^n$ term are formally determined 
and that there are no $L^2$ $\Delta_g$-harmonic functions.
\stopthm

We remark that on a general asymptotically hyperbolic manifold, the same
argument produces a unique solution of $\Delta_g U = n$ with leading
asymptotics $\log x +O(x)$.

If $n$ is odd, then $B|_M\in \CI(M)$ is determined by the choice of
Poincar\'e metric $g$ and representative metric $h$, and the
construction above shows that $B|_M=-\frac{d}{ds}S(s)1|_{s=n}$.  
We define an analogue of the $Q$-curvature for $n$ odd by
$$
Q=k_n B|_M=-k_n \frac{d}{ds}S(s)1|_{s=n}\,, \qquad k_n = 2^n
\frac{\Gamma(n/2)}{\Gamma(-n/2)}\,.  
$$
This quantity is globally determined and depends in general on the
extension of the formal Poincar\'e metric to a metric on $X$.  
It is tempting to try to
normalize this extension by requiring that it be exactly Einstein,
especially in light of the significant recent progress on this existence
question (\cite{A2}).   
However, an Einstein extension might not exist and might not in general be
unique.

For the odd dimensional version of $Q$, we have the following analogue of
\eqref{Qtrans}: 

\begin{proposition}\label{nodd}
If $n$ is odd, then under a conformal change ${\hat h} = e^{2\Up}h$, $Q$
satisfies the transformation law:
\begin{equation}\label{Qtransodd}
e^{n\Up}{\widehat Q} = Q +k_nS(n)\Up\,.
\end{equation}
\end{proposition}

\begin{proof}
As in the proof of \eqref{Qtrans} given in \S 3, the difference
${\hat U} - U$ is a solution of  
$\Delta_g({\hat U}-U)=0$ with $({\hat U} - U)|_M = \Up$.  By the definition
of the scattering matrix, the $x^n$ coefficient in the expansion of 
${\hat U} - U$ is $S(n)\Up$, which gives 
$e^{n\Up}\hat{B}|_M - B|_M = S(n)\Up$ as desired.
\end{proof}

In comparing \eqref{Qtransodd} with \eqref{Qtrans}, note that 
$k_nS(n)$ is a pseudodifferential operator on $M$ of order $n$, whose
principal symbol, by \eqref{symbol}, is equal to that of $\Delta_h^{n/2}$. 

Since $S(n)$ is self-adjoint and $S(n)1=0$, from \eqref{Qtransodd} it
follows that 
$\int_M Q\, dv_h$ is independent of the choice of $h$ in the conformal
class.  For the renormalized volume we have the following analogue of
Theorem~\ref{L}.

\begin{theorem}\label{V}
If $n$ is odd and $g$ is a Poincar\'e metric with conformal infinity $[h]$,
then the renormalized volume $V$ in \eqref{volexp} is given by
\begin{equation}\label{Vform}
k_nV= \int_M Q \,dv_h\,.
\end{equation}
\end{theorem}
\begin{proof}
Green's identity gives
$$
n\vol_g \{x>\ep\}= \int_{x>\ep}\Delta_g U \;dv_g = \ep^{1-n} \int_{x=\ep} 
\pa_xU\;dv_{h_{\ep}}\,.
$$
Evaluation of the constant term in the expansion after substituting 
\eqref{Uasym} and \eqref{volhx} gives $V=\int_MB|_M \, dv_h$ as desired. 
\end{proof}

\end{document}